
\input amstex
\documentstyle{amsppt}

\input label.def
\input degt.def

\input epsf
\def\picture#1{\epsffile{#1-bb.eps}}

\def\ie.{\emph{i.e\.}}
\def\eg.{\emph{e.g\.}}
\def\cf.{\emph{cf\.}}

\def\etc{\emph{etc}}

{\catcode`\@11
\gdef\proclaimfont@{\sl}}

\newtheorem{Remark}\Remark\thm\endAmSdef
\conjecture\thm\endproclaim
\problem\thm\endAmSdef

\def\dash{\item"\hfill--\hfill"}
\def\Dashes{\widestnumber\item{--}\roster}
\def\endDashes{\endroster}

\let\BAR\B

\loadbold
\def\bA{\bold A}
\def\bD{\bold D}
\def\bE{\bold E}
\def\bJ{\bold J}
\def\bU{\bold U}

\def\tA{\tilde\bA}
\def\tD{\tilde\bD}

\let\splus\oplus
\let\GAPcolon\rtimes

\let\Gb\beta
\let\Gd\delta
\let\Gs\sigma
\let\Go\omega
\def\CH{\Cal H}
\def\CC{\Cal C}
\def\CL{\Cal L}
\def\CM{\Cal M}
\def\CJ{\Cal J}
\def\CO{\Cal O}
\def\tS{\tilde S}
\def\tX{\tilde X}
\def\go{\frak o}
\def\LL{\Bbb L}
\def\FF{\Bbb F}
\let\tree\Xi

\def\Cp#1{\Bbb P^{#1}}
\def\Rp#1{\Cp{#1}_\R}

\def\Pic{\operatorname{Pic}}
\def\Aut{\operatorname{Aut}}
\def\discr{\operatorname{discr}}
\def\Sk{\operatorname{Sk}}

\let\<\langle
\let\>\rangle
\def\1{^{-1}}

\def\PGL#1{\operatorname{\text{\sl PGL}}_{#1}}
\def\MW#1{\operatorname{\text{\sl MW}}(#1)}
\def\O{\operatorname{\text{\sl O}}}
\def\SL{\operatorname{\text{\sl SL}}}
\def\PL{{\sl PL}}

\def\=#1={\sim_{\text{\rm#1}}}

\def\BG#1{\Bbb B_{#1}} 
\def\SG#1{\Bbb S_{#1}} 
\def\DG#1{\Bbb D_{#1}} 
\def\CG#1{\ZZ_{#1}}    
\def\QG{Q_8}           

\def\inserthyphen{\ifcat\next a-\fi\ignorespaces}
\let\BLACK\bullet
\let\WHITE\circ
\def\CROSS{\vcenter{\hbox{$\scriptstyle\mathord\times$}}}
\let\STAR*
\def\pblack-{$\BLACK$\futurelet\next\inserthyphen}
\def\pwhite-{$\WHITE$\futurelet\next\inserthyphen}
\def\pcross-{$\CROSS$\futurelet\next\inserthyphen}
\def\pstar-{$\STAR$\futurelet\next\inserthyphen}
\def\black{\protect\pblack}
\def\white{\protect\pwhite}
\def\cross{\protect\pcross}

\def\base{S^2}

\topmatter

\author
Alex Degtyarev
\endauthor

\title
Topology of plane algebraic curves:\\
the algebraic approach
\endtitle

\address
Department of Mathematics,
Bilkent University,
06800 Ankara, Turkey
\endaddress

\email
degt\@fen.bilkent.edu.tr
\endemail


\abstract
We discuss the principle tools and results and state a few open
problems concerning the classification and
topology of plane sextics and trigonal curves
in ruled surfaces.
\endabstract



\endtopmatter

\document

\section{Introduction}

\subsection{Contents}
The purpose of this survey is to discuss the algebro-geometric
aspects
of plane algebraic curves of small degree or
`codegree'
and state a few open problems concerning
topology of these curves.
More precisely, we will concentrate on two classes
of curves: sextics in the projective plane
and
trigonal curves in ruled surfaces, the latter being closely
related to plane curves
with a singular point of multiplicity
$(\text{degree}-3)$.

In the last $25$ years,
topology of singular plane algebraic
curves has been an area of active research;
the modern state of affairs and important open problems are
thoroughly
presented in recent surveys~\cite{Artal.survey},
\cite{Libgober.survey.1}, \cite{Libgober.survey.2},
and~\cite{Oka.survey}. Analyzing recent achievements, one
cannot help noticing that, with
relatively
few exceptions,
algebraic objects are studied by purely topological means, making
very little use of the analytic structure. These methods, largely
based on the braid monodromy, apply equally well to
pseudo-holomorphic curves and to the so called \emph{Hurwitz
curves}, see~\cite{Kulikov}, which roughly
are smooth surfaces mapped to~$\Cp2$
so that they behave
like algebraic curves with respect to one chosen pencil
of lines. (Remarkably, a similar phenomenon was observed in
topology of real algebraic curves, leading Rokhlin
and Viro to the concept
of \emph{flexible curves}. This idea was later developed by
Orevkov, who initiated the study of real pseudo-holomorphic
curves.) Although well justified in general,
this approach seems to fail when curves of small degree are
concerned and more precision is required. For this reason, I chose
to restrict this survey to plane sextics and to the more algebraic
methods that have recently been used in their study.
As part of these methods, trigonal curves appear; since they seem
to be of interest in their own right,
I dedicated a separate section
to their theory and related problems.

Although
this paper is of a purely expository nature,
I do
announce several new results, among them
being Theorem~\ref{th.es.tree},
computing the transcendental lattices of
a series of elliptic surfaces,
and updated lists of fundamental groups of sextics, see comments
to Problems~\ref{prbl.pi1} and~\ref{prbl.tc.pi1}.

\subsection{Acknowledgements}
I would like to thank J.~I.~Cogolludo and E.~Hironaka for
organizing a remarkable conference and encouraging me to write
this survey.
I am also grateful to Vik\.~S.~Kulikov, S.~Orevkov,
I.~Shimada, E.~Shustin,
and A.~Vainshtein for expanding my horizons on various related
subjects, as well as to
the speakers and participants of the
Lib{\it60\/}ber conference for the enlightening talks and
fruitful discussions.

\section{Preliminaries}\label{S.classical}

Throughout the paper,
all varieties are either complex analytic (when speaking about the
moduli of $K3$-surfaces)
or algebraic
over~$\C$. Unless stated otherwise, all curves considered are reduced.

\subsection{Equivalence relations}\label{s.equivalence}
The principal subject of the paper is the topology of a reduced
plane algebraic curve $C\subset\Cp2$ or, more generally, a curve
$C$ in a fixed algebraic surface.
Typically, plane curves are considered up to one of the following major
equivalence relations:
\roster
\item\local{equiv.def}
$\=def=$ -- \emph{\rom(equisingular\rom) deformation equivalence}
(or \emph{rigid isotopy}):
$C_1\=def=C_2$ if $C_1$ and~$C_2$ can be connected by a path in an
equisingular stratum of the space~$\CC_d$ of plane algebraic curves
of a given degree~$d$;
\item\local{equiv.PL}
$\=PL=$ -- \emph{\PL-equivalence}:
$C_1\=PL=C_2$ if the pairs $(\Cp2,\C_1)$ and $(\Cp2,C_2)$ are
\PL-homeomorphic;
\item\local{equiv.top}
$\=top=$ -- \emph{topological equivalence}:
$C_1\=top=C_2$ if the complements
$\Cp2\sminus C_1$ and $\Cp2\sminus C_2$ are homeomorphic;
\item\local{equiv.cfg}
$\=cfg=$ -- \emph{combinatorial equivalence}: $C_1\=cfg=C_2$ if
$C_1$ and $C_2$ have the same combinatorial type of singularities
(or \emph{configuration type}); formally, this means that
pairs $(T_1,C_1)$
and $(T_2,C_2)$ are \PL-homeomorphic, where
$T_i$ is a regular neighborhood of~$C_i$ in~$\Cp2$,
$i=1,2$.
\endroster
Clearly, one has
\loccit{equiv.def}$\implies$\loccit{equiv.PL}$\implies$\loccit{equiv.top}
and~\loccit{equiv.cfg}.

\Remark
In order to avoid moduli in the case of
non-simple singular points, we use \PL-homeomorphisms rather than
diffeomorphisms in Items~\itemref{s.equivalence}{equiv.PL}
and~\itemref{s.equivalence}{equiv.cfg}. For the same reason,
`equisingular' in Item~\itemref{s.equivalence}{equiv.def}
is understood in the
\PL-category. (Alternatively,
two singular points are considered equivalent if they have
combinatorially isomorphic resolutions.)
\endRemark

\Remark
The configuration type, see Item~\itemref{s.equivalence}{equiv.cfg},
of an irreducible curve is determined by
its degree and the set of \PL-types of its singularities. In general,
one should take into account the degrees of the components of the
curve and the distribution of the branches at the singular points
among the components. For a precise combinatorial definition,
see~\cite{ArtalTokunaga}.
\endRemark

\Remark
In Item~\itemref{s.equivalence}{equiv.top}, one could as
well consider homeomorphisms of the pairs $(\Cp2,C_i)$, \cf.
Item~\itemref{s.equivalence}{equiv.PL}.
However,
most invariants currently
used to distinguish curves take into account
the complement $\Cp2\sminus C_i$ only.
\endRemark

\subsection{Classical problems and results}\label{s.classical}
The two principal problems of topology of plane algebraic curves
are the classification of curves (of a given degree, with a given
set of singularities, \etc\.) up to one of the equivalence
relations discussed in Subsection~\ref{s.equivalence} and the
study of the topological type of the pair $(\Cp2,C)$ or the
complement $\Cp2\sminus C$. In the latter case, of special
interest is the fundamental group $\pi_1(\Cp2\sminus C)$, which is
commonly referred to as the \emph{fundamental group of~$C$}; this
group controls the existence of algebraic surfaces ramified
at~$C$, reducing in a sense the study of surfaces to that of
curves. For the particular groups appearing in this paper, we use
the following notation:
\Dashes
\dash
$\CG{n}$ is the cyclic group of order~$n$;
\dash
$\DG{2n}$ is the dihedral group of order~$2n$;
\dash
$\SG{n}$ is the symmetric group of degree~$n$;
\dash
$\BG{n}$ is the braid group on $n$ strands, and
$\Gs_1,\ldots,\Gs_{n-1}$ is its Artin basis;
\dash
$\QG:=\{\pm1,\pm i,\pm j,\pm k\}\subset\Bbb H$ is the
\emph{quaternion group} of order~$8$.
\endDashes
In general, the group $\pi_1(\Cp2\sminus C)$ is very difficult to
compute.
Following Zariski,
Libgober~\cite{Libgober1}
suggested an intermediate invariant, the \emph{Alexander
polynomial} $\Delta_C(t)$: if $C$ is irreducible, $\deg C=d$, one
can consider a resolution of singularities $\tX_d$ of
the $d$-fold cyclic covering $X_d\to\Cp2$ ramified
at~$C$ and define $\Delta_C(t)$ as the characteristic polynomial
of the deck translation action on $H_1(\tX_d;\C)$.

As mentioned in the Introduction, most tools used in the subject
are of a purely topological nature and apply to arbitrary Hurwitz
curves. The two truly algebraic results that I am aware of are
Nori's theorem~\cite{Nori} (which is the best known generalization
of Zariski's conjecture on nodal curves) and a formula computing
the Alexander polynomial, see Theorem~\ref{th.Delta} below.

For Nori's theorem, introduce the invariant $n(S)$ of a (type of)
singularity~$S$ as follows: $n(\bA_1)=2$, and for any other type,
$n(S)$ is the decrement of the self-intersection of the curve when
$S$ is resolved to a divisor with normal intersections.
Then the following statement holds, see~\cite{Nori}.

\theorem\label{th.Nori}
Given an irreducible plane curve~$C$, if
$C^2>\sum n(S)$, the summation running over all singular
points~$S$ of~$C$,
then the group
$\pi_1(\Cp2\sminus C)$ is abelian.
\endtheorem

For the Alexander polynomial, given a singular germ~$S$ and a
rational $r\in\Q$, one can define the so called \emph{ideal of
quasiadjunction}
$\CJ_{S,r}\subset\CO_S$, depending only on~$r$ and the
(analytic) type of~$S$. (For the precise definition and
computation of $\CJ_{S,r}$, see the citations
at the end of this paragraph.)
Now, given an algebraic curve
$C\subset\Cp2$ of degree~$d$ and an integer
$i\in\{0,\ldots,d-1\}$, define the \emph{quasiadjunction sheaf
of ideals}
$\CJ_i\subset\CO$ as follows:
the stalk $\CJ_i|_S$ equals $\CJ_{S,i/d}$
at each singular point~$S$ of~$C$, and $\CJ_i|_P=\CO_P$ at any
other point~$P$. Finally, define the linear system~$\CL_i$ as the
space of sections of $\CJ_i\otimes\CO(i-3)$.
Next theorem, generalizing Zariski's
computation~\cite{Zariski.group} for the six cuspidal sextic,
appeared under various disguises in~\cite{poly}, \cite{Esnault},
\cite{Libgober.th}, \cite{Loeser}, and~\cite{Oka.poly}.

\theorem\label{th.Delta}
In the notation above, one has
$$
\Delta_C(t)=\prod_{i=0}^r[(t-\xi^i)(t-\xi^{-i})]^{\dim H^1(\Cp2;\CL_i)},
$$
where $\xi_i=\exp(2\pi\sqrt{-1}/d)$ is the $d$-th primitive root
of unity.
\endtheorem

In classical terms, the number $\dim H^1(\Cp2;\CL_i)$ is called
the \emph{superabundance} of~$\CL_i$; it is equal to the
difference between the actual dimension of~$\CL_i$ and its virtual
(or expected) dimension. (Note that the
expected dimension depends on the
degree and the set
of singularities of~$C$ only.)
As an interesting consequence,
the dimensions of certain linear systems turn out to be
topological invariants. A similar phenomenon is discussed in the
comments to Problems~\ref{prbl.splitting}
and~\ref{prbl.splitting.d} below.

Theorem~\ref{th.Delta} admits generalizations to
curves/divisors in arbitrary surfaces; there are versions of this
theorem for the multivariable Alexander polynomials and
characteristic varieties. For further references,
see~\cite{Libgober.survey.1} and~\cite{Oka.survey}.

\section{Plane sextics: tools and results\label{S.sextics}}

Simple sextics are related to $K3$-surfaces; we briefly
outline this
relation and its implications to topology and geometry of sextics
in Subsections~\ref{s.simple} and~\ref{s.K3}.
Non-\allowbreak{}simple sextics
are dealt with in Subsection~\ref{s.nonsimple}.


\subsection{Classes of sextics}\label{s.6.intro}
First, we divide reduced plane sextics into \emph{simple} and
\emph{non-simple}, the former having simple (\ie., $\bA\bD\bE$
type) singularities only and the latter having at least one
non-simple singular point. Simple sextics are related to
$K3$-surfaces, whereas non-simple ones are related to rational
and, in few cases, irrational ruled surfaces. As a consequence,
simple sextics turn out to be much more complicated;
they are the subject of most problems below.

Another important class of sextics, first discovered by
Zariski~\cite{Zariski.group}, is formed by the so called sextics
of torus type.

\definition\label{def.torus}
A sextic $C\subset\Cp2$ is said to be of \emph{torus type} if one
of the following two equivalent conditions holds:
\roster
\item\local{torus.pq}
$C$ is given by an equation of the form $f_2^3+f_3^2=0$, where $f_2$
and~$f_3$ are homogeneous polynomials of degree~$2$ and~$3$,
respectively;
\item\local{torus.locus}
$C$ is the ramification locus (assuming that it is reduced) of a
projection to~$\Cp2$ of a cubic surface in~$\Cp3$.
\endroster
If $C$ is irreducible, conditions~\loccit{torus.pq},
\loccit{torus.locus} are also equivalent to any of the following
three conditions (see~\cite{degt.Oka} and~\cite{degt.Oka2}):
\roster
\item[3]\local{torus.B3}
the group $\pi_1(\Cp2\sminus C)$ factors to the
reduced braid group $\BG3/(\Gs_1\Gs_2)^3$;
\item\local{torus.S3}
the group $\pi_1(\Cp2\sminus C)$ factors to the
symmetric group $\SG3=\DG6$;
\item\local{torus.Delta}
the Alexander polynomial $\Delta_C(t)$ is nontrivial,
$\Delta_C(t)\not\equiv1$.
\endroster
\enddefinition

The fact that condition~\loccit{torus.Delta} is equivalent
to~\loccit{torus.pq} was first conjectured by
Oka~\cite{Oka.conjecture}.

Remarkably, conditions~\loccit{torus.pq} and~\loccit{torus.locus}
above remain invariant under equisingular deformations. The
simplest sextic of torus type, discovered in~\cite{Zariski.group},
has six cusps; it is the ramification locus of a generic
projection of a nonsingular cubic, and its fundamental group is
$\BG3/(\Gs_1\Gs_2)^3\cong\CG2*\CG3$. This curve and its properties
can be generalized in at least two directions. First, one can
consider the so called \emph{curves of $(p,q)$-torus type}, \ie.,
curves of degree~$pq$ given by an equation
$f_q^p+f_p^q=0$,
see \cite{Oka.deg10},
\cite{Oka.Zariski}, \cite{Oka.torus.1},
and~\cite{Kulikov.torus} for a further generalization.
Second, one can study the ramification locus of a generic
projection to~$\Cp2$ of a hypersurface in~$\Cp3$
(see~\cite{Moishezon}) or, more generally, a smooth surface
$X\subset\Cp{N}$, see
~\cite{Teicher.survey} for
further references.
Both approaches provide interesting examples of plane curves.

\subsection{Sextics with a non-simple singular point}\label{s.nonsimple}
Non-simple sextics turn out to be `simple',
and most classical problems
related to such curves can be solved by
fairly elementary methods.
The principal results and tools used are outlined here.

\theorem\label{th.6.nonsimple}
The equisingular deformation type of an irreducible sextic
with a non-simple singular point is determined by its
set of singularities\rom; in other words, the relations $\=def=$
and $\=cfg=$ restricted to such sextics are equivalent.
\endtheorem

A complete proof of this theorem, although rather straightforward,
has never been published. Various partial statements, as well as a
complete
description of the combinatorial types realized by irreducible
sextics, are scattered
through~\cite{quintics}, \cite{degt.Oka}, and~\cite{degt.Oka2}. An
independent treatment of sextics of torus type, including
non-simple ones, is found
in~\cite{OkaPho.moduli}, \cite{OkaPho}, and~\cite{Pho}.

The fundamental groups of non-simple sextics
are also all known (in the sense that their presentations can
easily be written down); for irreducible curves, the groups are
listed in~\cite{degt.Oka} and~\cite{degt.Oka2}. Apart from the
$21$ deformation families of
sextics of torus type
(all but two groups being $\BG3/(\Gs_1\Gs_2)^3$),
there are five families with
$\pi_1(\Cp2\sminus C)=\DG{10}\times\CG3$; all other groups are
abelian.

A non-simple singular point has multiplicity at least three, and
the projection from this point restricts to a degree at most~$3$
map $C\to\Cp1$. The braid monodromy of this projection
is easily computable and
gives one a fairly good understanding of the
topology of the pair $(\Cp2,C)$. The most involved (and relatively
new) is the case of
multiplicity three, \ie., a singular point adjacent to~$\bJ_{10}$ in
Arnol$'$d's notation (simple tangency of three smooth branches).
The corresponding sextics are so called
\emph{trigonal curves}; they can be treated as described in
Section~\ref{S.trigonal} below.

These techniques should work for reducible
non-simple sextics as well. However, in view of the number of
details to be taken into account, I would refrain from a precise
statement and pose it as a problem. (For other applications, see
also comments to Problem~\ref{prbl.degenerations}.)

\problem\label{prbl.non-simple}
Classify reducible non-simple sextics. Are the relations $\=def=$
and $\=cfg=$
equivalent for such sextics?
\endproblem

\subsection{Classification of simple sextics}\label{s.simple}
Recall that a \emph{lattice} is a
free abelian
group~$L$ of finite rank
supplied with a symmetric bilinear form
$L\otimes L\to\Z$ (which is usually referred to as \emph{product}
and denoted by $x\otimes y\mapsto x\cdot y$ and
$x\otimes x\mapsto x^2$). A \emph{root} in an even
lattice is a vector
of square~$(-2)$; a \emph{root system} is a negative definite
lattice generated by its roots.

It is convenient to identify a set of simple singularities with a
root system, replacing each singular point with
the irreducible root system of the same name ($\bA_p$, $\bD_q$,
$\bE_6$, $\bE_7$, or $\bE_8$) and taking a direct sum
over all singular points in the set.
A \emph{marking} of a set of simple
singularities~$\Sigma$ is a choice of a fundamental system of
roots $\Gs\subset\Sigma$.
(Recall that $\Gs$ can be defined as the
set of walls of a fixed Weyl chamber of~$\Sigma$; the elements
of~$\Gs$ form a basis for~$\Sigma$, so that $\Sigma=\Z\Gs$.)

\definition\label{def.config}
A \emph{configuration}
(a set of \emph{lattice data} in~\cite{Shimada.MMJ})
extending a given marked
set of simple singularities $(\Sigma,\Gs)$ is a finite
index lattice extension $\tS\supset S:=\Sigma\oplus\Z h$, $h^2=2$,
satisfying the following conditions:
\roster
\item\local{config.no.roots}
each root $r\in\tS\cap(\Sigma\otimes\Q)$ belongs to~$\Sigma$;
\item\local{config.not.elliptic}
there is no root $r\in\Sigma$ such that $\frac12(r+h)\in\tS$.
\endroster
An \emph{isomorphism} of two configurations
$\tS'\supset\Z\Gs'\oplus\Z h'$ and
$\tS''\supset\Z\Gs''\oplus\Z h''$ is a lattice isomorphism
$\tS'\to\tS''$ taking~$h'$ to~$h''$ and $\Gs'$ to~$\Gs''$.
\enddefinition

From now on, we fix a lattice $\LL\cong2\bE_8\oplus3\bU$, where
$\bU$ is the \emph{hyperbolic plane},
\ie., the lattice spanned by two
generators $u$, $v$ with $u^2=v^2=0$ and $u\cdot v=1$.

\definition\label{def.CH}
An \emph{abstract homological type} extending a given
configuration $\tS\supset S:=\Sigma\oplus\Z h$
is an isometry $\CH\:\tS\to\LL$ taking~$\tS$ to a primitive
sublattice of~$\LL$.
An \emph{orientation} of~$\CH$
is an orientation $\go$ of maximal positive definite
subspaces in the orthogonal complement $\CH(\tS)^\perp$.
An \emph{isomorphism} of two
(oriented)
abstract homological types
$\CH'\:\tS'\to\LL$ and $\CH''\:\tS''\to\LL$ is an automorphism
of~$\LL$ inducing an isomorphism $\tS'\to\tS''$ of configurations
(and taking the orientation of~$\CH'$ to the
orientation of~$\CH''$).
\enddefinition

\Remark\label{rem.Weyl}
Any two fundamental systems of roots
in~$\Sigma$ can be interchanged by a sequence of reflections
of~$\Sigma$, which would extend trivially to any larger lattice
containing~$\Sigma$. Hence, for the simple enumeration tasks
(\cf. Theorem~\ref{th.6.enum} below), the marking can be ignored.
It is, however, important in some more subtle problems, \eg., in
the study of automorphisms of an abstract homological type.
\endRemark

Let $C\subset\Cp2$ be a simple sextic.
Consider the double covering $X\to\Cp2$
ramified at~$C$ and its minimal resolution
$\tX=\tX_C\to X$. It is a $K3$-surface.
Introduce the following objects:
\Dashes
\dash
$\Gs_C\subset H_2(\tX)$,
the set of classes of the exceptional
divisors contracted in~$X$;
\dash
$\Sigma_C\subset H_2(\tX)$,
the sublattice spanned by~$\Gs_C$;
\dash
$h_C\in H_2(\tX)$,
the pull-back of the class $[\Cp1]\in H_2(\Cp2)$;
one has $h_C^2=2$;
\dash
$\tS_C\subset H_2(\tX)$,
the primitive hull of $S_C:=\Sigma_C\oplus\Z h_C$
in $H_2(\tX)$;
\dash
$T_C=S_C^\perp$, the stable
\emph{transcendental lattice} of~$\tX$;
\dash
$\Go_C\in T_C\otimes\C$, the class of a holomorphic form
on~$\tX$ (the \emph{period} of~$\tX$);
\dash
$\go_C$, the orientation of the plane
$\Go_C^\R:=\R\Re\Go_C\oplus\R\Im\Go_C\subset T_C\otimes\R$.
\endDashes
Then $\Sigma_C$ is isomorphic to the set of singularities of~$C$
(regarded as a lattice), $\Gs_C$ is a marking of~$\Sigma_C$,
and the extension
$\tS\supset S$ is a configuration.
Hence, fixing an isomorphism $\phi\:H_2(\tX)\to\LL$
(called a
\emph{marking} of~$\tX$) and combining it with the inclusion
$\tS_C\hookrightarrow H_2(\tX)$, one obtains an abstract
homological type $\CH_C\:\tS_C\to\LL$, called the
\emph{homological type} of~$C$. The orientation $\go_C$ extends to
an orientation of~$\CH_C$. The pair $(\CH_C,\go_C)$ is well
defined up to the choice of a marking~$\phi$,
\ie., up to isomorphism.

Next statement is found in~\cite{JAG}; the
existence part was first proved
in~\cite{Urabe} and then exploited in~\cite{Yang}.

\theorem\label{th.6.enum}
The map $C\mapsto(\CH_C,\go_C)$ establishes a
bijection
between the set of equisingular deformation classes
of simple sextics and the set of isomorphism classes of oriented
abstract homological types.
\endtheorem

The proof of Theorem~\ref{th.6.enum}, as well as of most other
statements related to simple sextics, is based on the so called
global Torelli theorem and surjectivity of the period map for
$K3$-surfaces, see \eg.~\cite{BPV}. A convenient restatement of
these two facts and a description of a fine moduli space
of marked polarized $K3$-surfaces is found in~\cite{Beauville}. As
another consequence, one obtains a description of the
equisingular moduli spaces.
Denote by $\CC_6(\CH,\go)$ the equisingular stratum
of the space of sextics defined by
an oriented abstract homological type $(\CH,\go)$, and let
$\CM_6(\CH,\go)=\CC_6(\CH,\go)/\!\PGL3$ be the corresponding moduli
space.
Let $T=\CH(\tS)^\perp$, and let~$\Omega_{\CH}$ be the
projectivization of the cone
$\{\Go\in T\otimes\C\,|\,\Go^2=0,\ \Go\cdot\bar\Go>0\}$. This
space has two complex conjugate connected components, which differ
by the orientation of the positive definite $2$-subspace
$\Go^\R:=\R\Re\Go\oplus\R\Im\Go\subset T\otimes\R$; denote
by~$\Omega_{\CH}^+$ the one that matches~$\go$. Then the following
statement holds.

\theorem\label{th.moduli}
The equisingular moduli space $\CM_6(\CH,\go)$
can be identified with
a Zariski open subset of the
quotient
$\Omega_{\CH}^+/\!\Aut(\CH,\go)$.
In particular, $\CM_6(\CH,\go)$ is irreducible and one has
$\dim\CM_6(\CH,\go)=19-\rank\Sigma$.
\endtheorem

Note that, for a sextic~$C$, the rank $\rank\Sigma_C$ equals the
\emph{total Milnor number} $\mu(C)$; roughly, it is the number of
algebraic conditions selecting curves with the prescribed set of
singularities. It follows that for simple sextics these conditions
are always independent (\cf. Theorem~\ref{th.perturbations}
below).

A sextic~$C$ with $\mu(C)=19$
is called \emph{maximizing}. The moduli space of
maximizing sextics is discrete, \ie., two maximizing sextics are
deformation equivalent if and only if they are related by a
projective transformation.
Any maximizing sextic is defined over an algebraic
number field (see comments to Problem~\ref{prbl.Galois} below).

\subsection{Further applications of theory of
$K3$-surfaces}\label{s.K3}
As another consequence of the description of the moduli space of
$K3$-surfaces, one obtains the following statement, concerning the
generic (within an equisingular family) Picard group.

\theorem\label{th.Picard}
There is a dense Zariski open subset
$\CC'\subset\CC_6(\CH,\go)$ such
that, for a sextic $C\in\CC'$,
the Picard group $\Pic\tX_C$ equals $\tS_C$.
\endtheorem

In~\cite{Shimada.MMJ}, a sextic~$C$ with $\Pic\tX_C=\tS_C$ is
called \emph{lattice generic}. In fact, for a lattice generic
sextic~$C$ one can even describe the K\"ahler cone of~$\tX_C$
(hence the semigroup of numerically effective divisors) and the
set of (the classes of) irreducible rational curves on~$\tX_C$,
see \eg.~\cite{BPV}; this description makes use of the fundamental
polyhedron of the group generated by reflections of the hyperbolic
lattice~$\tS$, see~\cite{Vinberg}. The key r\^ole is played
by the Riemann--Roch theorem, which, in the case of $K3$-surfaces,
takes an especially simple form and asserts that, for any
divisor~$D$ with $D^2\ge-2$, either $D$ or $-D$ is effective.

Theorem~\ref{th.Picard} means that the configuration
$\tS_C\supset S_C$ captures the essential algebro-geometric
information invariant under equisingular
deformations. This fact lead Shimada~\cite{Shimada.MMJ} to the
introduction of another equivalence relation, the so called
lattice equivalence.

\definition\label{def.lattice}
Two simple sextics $C_1$, $C_2$ are \emph{lattice equivalent},
$C_1\=lat=C_2$, if their configurations
$\tS_{C_1}\supset S_{C_1}$ and $\tS_{C_2}\supset S_{C_2}$ are
isomorphic.
\enddefinition

This new relation is strictly between $\=def=$ and $\=cfg=$
(see~\cite{Yang}) and at present it seems somewhat easier to
handle. The following geometric properties of simple sextics are
known to be lattice invariant: the existence of certain dihedral
covering and, in particular, being of torus type
(see~\cite{degt.Oka} and~\cite{Tokunaga.new}),
the existence of stable symmetries (in the irreducible case,
see~\cite{symmetric} and Problem~\ref{prbl.symmetries}
below), the existence of $Z$-splitting lines,
conics, and cubics (see~\cite{Shimada.MMJ}
and comments to
Problems~\ref{prbl.splitting} and~\ref{prbl.splitting.d}
below).
Most of these
statements (except the existence of dihedral coverings) are based
on the Riemann--Roch theorem
and do not extend to higher degrees.

Another remarkable result is a clear picture of the adjacencies of
the equisingular strata of the space of simple sextics. According
to~\cite{Looijenga}, the isomorphism classes of perturbations of a
simple singularity are enumerated by the induced subgraphs of its
Dynkin diagram (up to a certain equivalence). Given an abstract
homological type $\CH\:\tS\to\LL$,
$\tS\supset S\supset\Sigma\supset\Gs$,
define its \emph{combined Dynkin diagram}~$D_{\CH}$ as the
Dynkin diagram of the fundamental system of roots~$\Gs$. Any induced
subgraph $D'\subset D_{\CH}$ gives rise to a fundamental system of
roots $\Gs'\subset\Gs$ and hence to a configuration
$\tS':=\tS\cap(S'\otimes\Q)\supset S'$, where
$S'=\Z\Gs'\oplus\Z h$. Restricting $\CH$ to~$\tS'$,
one obtains a new abstract homological type~$\CH'$, and any
orientation~$\go$ of~$\CH$ extends to an orientation~$\go'$
of~$\CH'$ in the obvious way.
We call the pair $(\CH',\go')$ the \emph{restriction} of
$(\CH,\go)$ to the subgraph $D'\subset D_{\CH}$.
The following statement is contained
in~\cite{degt.8a2} (see also~\cite{Shimada.MMJ} for a more
formal proof).

\theorem\label{th.perturbations}
A stratum $\CC_6(\CH,\go)$ is in the closure
of $\CC_6(\CH',\go')$
if and only if $(\CH',\go')$ is isomorphic to the
restriction of $(\CH,\go)$ to an induced subgraph
$D'\subset D_{\CH}$.
\endtheorem

In other words, the singular points of a simple sextic can be
perturbed arbitrarily and independently. This fact is in sharp
contrast to the case of higher degrees. For example, for each of
the following curves
\Dashes
\dash
the Pappus configuration of nine lines,
\dash
union of three cubics passing through nine common points,
\dash
union of a three cuspidal quartic and the tangents at its three
cusps,
\endDashes
one cannot perturb a triple point to three nodes while keeping the
other singular points intact. (These examples were communicated to
me by E.~Shustin; note that all curves have simple singularities
only.) One can construct a great deal of other
examples in degrees~$7$ and~$8$ by combining a simple sextic with
a $Z$-splitting line or conic (see~\cite{Shimada.MMJ} and comments
to Problems~\ref{prbl.splitting} and~\ref{prbl.splitting.d}
below).
At present, no example if irreducible curve with restricted
perturbation is known, but there is no reason to believe that such
examples do not exist. Some general conditions sufficient for the
singular points of a plane curve to be perturbed independently are
found in~\cite{Shustin}; however, these conditions are too weak to
cover all simple sextics.

It is not difficult to prove an analogue of
Theorem~\ref{th.perturbations} for non-simple sextics, {\em provided
that the non-simple points are kept non-simple and their
multiplicity does not change} (so that the projection mentioned in
Subsection~\ref{s.nonsimple} is preserved during the
perturbation; this projection can be used to describe the
degenerations). This observation gives rise to the following
problem.

\problem\label{prbl.degenerations}
Describe the degenerations of simple sextics to non-simple ones,
as well as the degenerations increasing the multiplicity of a
non-simple point.
\endproblem

A partial answer to this question is given in~\cite{OkaPho.moduli}
and~\cite{OkaPho}, where degenerations are used in the study
of sextics of torus type.

\section{Simple sextics: problems\label{S.simple}}

One can anticipate that the two principal problems concerning
sextics, namely
their classification and the computation of their fundamental
groups, will be completely solved within a few years.
However, there still
remain more subtle geometric questions, see
Subsection~\ref{s.geometry}, and
attempting
to understand the extent to which
the reach experimental material gathered for sextics reflects
properties of algebraic curves in general
would be of utmost interest.

\subsection{The classification}\label{s.classification}
We start with a few classification questions.

\problem\label{prbl.classification}
Complete the deformation classification of simple sextics.
\endproblem

Theorem~\ref{th.6.enum} reduces this problem to a purely
arithmetical question, which can be solved using Nukulin's
theory~\cite{Nikulin} of discriminant forms. As a preliminary
step, Yang~\cite{Yang} has compiled a complete list of combinatorial
classes and Shimada~\cite{Shimada.MMJ} has compiled a complete list of
lattice equivalence classes; the latter contains about 11,500
items. Besides, Shimada~\cite{Shimada} has also listed all
deformation classes with the maximal total Milnor number $\mu=19$.
In the case $\mu\le18$, one expects to have very few
configurations extending to more than one oriented abstract
homological type. (At present, only one such example is known:
irreducible sextics of torus type with the set of singularities
$\bE_6\oplus\bA_{11}\oplus\bA_1$ form two complex conjugate
deformation families, see~\cite{Aysegul}.) However, in this case
one needs to deal with the automorphism
groups of indefinite lattices, which are not very well known.
One can hope that most such
lattices can be
handled using the results of~\cite{Miranda1} and~\cite{Miranda2},
which give a precise description of the cokernel of the natural
homomorphism $\O(L)\to\Aut\discr L$ for an indefinite lattice~$L$
of rank at least three.

It is expected that the final list will contain $11.5$ to $12$
thousand classes. Once the list is completed, next task
would be understanding the
result and deriving geometric consequences.
For example, as a rule, maximizing
sextics are much easier to construct explicitly,
and having computed
their topological invariants (such as the braid monodromy and
fundamental group, see comments to
Problem~\ref{prbl.pi1} below),
one can obtain the invariants of all perturbed
curves. Hence, the following question is of utmost importance.

\problem\label{prbl.degeneration}
Find a complete list of (irreducible) simple sextics that do
\emph{not} admit a degeneration to an (irreducible) maximizing
one.
\endproblem

This question can probably be answered using the information
already at hand, namely Shimada's list~\cite{Shimada.MMJ} of
lattice types.
One example is an irreducible sextic with the set of
singularities~$9\bA_2$: it does not admit any further degeneration
(in the class of simple sextics).

\subsection{The topology}\label{s.topology}
Next group of problems concerns the topology and homotopy
type (in particular, the fundamental
group) of the complement of a simple sextic.

\problem\label{prbl.Zariski}
Do all pairs of not deformation equivalent (irreducible)
sextics sharing the same configuration type form Zariski
pairs in the sense of Artal~\cite{Artal.def}, \ie., do they differ
topologically?
\endproblem

It is worth mentioning that, in the case of
simple singularities only, one can replace \PL-homeomorphisms with
diffeomorphisms in the definitions of
equivalence relations in Subsection~\ref{s.equivalence};
thus, it follows from~\cite{JAG} that
the relations $\=def=$ and $\=PL=$
(see~\iref{s.equivalence}{equiv.def}
and~\ditto{equiv.PL}, respectively) are equivalent when restricted
to simple sextics. Hence, the existence of a pair $C_1$, $C_2$
with $C_1\not\=def=C_2$ but $C_1\=top=C_2$ would reflect some
subtle phenomena of topology of $4$-manifolds.

Several dozens of candidates to be tried are readily found in Shimada's
list~\cite{Shimada} of maximizing sextics: one should consider the
pairs~$C_1$, $C_2$ that share the same configuration and
transcendental lattice~$T_{C_i}$ but are not complex conjugate.
Some of these pairs have been intensively studied,
see~\cite{Artal.braids}, \cite{Artal.Trends}, \cite{Artal.KT},
\cite{degt.e6}, \cite{dessin-e7}, \cite{dessin-e8},
\cite{dessin-e6}, \cite{EyralOka},
and the survey~\cite{Artal.survey} for further
references.
Most pairs are either known or expected to be
Galois conjugate, see Problem~\ref{prbl.Galois} below.
In many cases, the fundamental groups
$\pi_1(\Cp2\sminus C_i)$ have been computed; they are either equal
(when finite) are are not known to be distinct (irreducible
sextics of torus type with the sets of singularities
$2\bE_6\splus\bA_5\splus\bA_2$,
$\bE_6\splus\bA_5\splus\bA_4\splus2\bA_2$,
$\bE_6\splus\bA_8\splus2\bA_2\splus\bA_1$,
$\bE_6\splus2\bA_5\splus\bA_3$;
the groups must have isomorphic profinite completions, which makes
it difficult to distinguish them).

Within each pair~$C_1$, $C_2$ as above, the two curves differ by a
rather subtle lattice theoretic invariant, which is responsible
for one of the homomorphisms in the Mayer--Vietoris exact sequence
of the triad $(\tX_i\sminus C_i,N_i;\partial N_i)$, where $N_i$ is
a regular neighborhood of~$C_i$ in $\tX_i=\tX_{C_i}$, $i=1,2$.
As a first attempt, one can try to assume that
$\pi_1(\Cp2\sminus C_i)=\CG6$ and draw conclusions about the
homotopy types of the complements $\Cp2\sminus C_i$ or
the $6$-fold cyclic
coverings $Y_i\to\Cp2\sminus C_i$ by computing the homology
of~$Y_i$ (preferably as
$\pi_1(\Cp2\sminus C_i)$-modules).

\problem\label{prbl.pi1}
What are the fundamental groups of irreducible simple sextics?
\endproblem

There is a vast literature on the subject,
which calls for a separate survey.
The current
conjecture, suggested in~\cite{degt.Oka3}, is the following.

\conjecture\label{conj.finite}
The fundamental group of an irreducible sextic
that is not of torus type is finite.
\endconjecture

Note that the group of a sextic of torus type is never
finite, \cf. condition~\iref{def.torus}{torus.B3}.
Conjecture~\ref{conj.finite} replaces original Oka's
conjecture~\cite{Oka.conjecture} (its part concerning the
fundamental group) that was disproved in~\cite{degt.Oka}.
The conjecture has been verified for about 1,500 deformation
families not covered by Nori's theorem~\cite{Nori}. In particular,
it is known to be true for the following classes:
\Dashes
\dash
all non-simple sextics (\cite{degt.Oka},~\cite{degt.Oka2}),
\dash
all maximizing
sextics with an $\bE$~type singular point
(\cite{Artal.Trends}, \cite{dessin-e7},
\cite{dessin-e8},~\cite{dessin-e6}),
\dash
all but one
sextics whose group admits a dihedral quotient
(\cite{degt.Oka3}, \cite{degt-Oka}, ~\cite{Oka.D10}).
\endDashes
A number
of sporadic examples is contained in~\cite{Artal.braids},
\cite{Artal.Trends}, \cite{Artal.KT},
\cite{Oka.conjecture}, \cite{EyralOka1}, \cite{EyralOka2}, and a
huge number of other curves can be obtained from those already
mentioned by using Theorem~\ref{th.perturbations}. So far, very
few nonabelian fundamental
groups have been found,
their commutants being
$$
\CG5,\ \
\CG7,\ \
\SL(2,\CG9)\GAPcolon\CG5,\ \
((\CG2\times\QG)\GAPcolon\CG2)\GAPcolon\CG5,\ \
\SL(2,\FF_5),\ \
\SL(2,\FF_{19}).
$$
With one possible exception (the set of singularities
$3\bA_6\splus\bA_1$), only three soluble groups can appear:
$$
\CG3\times\DG{10},\ \
\CG3\times\DG{14},\ \
(((\CG2\times\QG)\GAPcolon\CG2)\GAPcolon\CG5)\GAPcolon\CG6.
$$
(Note that, being soluble, $\pi_1(\Cp2\sminus C)$ must admit a
dihedral quotient.)

The fundamental groups of most
(all but ten sets of singularities)
irreducible sextics of torus type
are also known, see~\cite{degt.2a8} for a `map' of the results
and~\cite{dessin-e6} for a recent update. With few exceptions,
they are minimal possible, \ie., $\BG3/(\Gs_1\Gs_2)^3$, \cf.
condition~\iref{def.torus}{torus.B3}.

In order to compute the group $\pi_1(\Cp2\sminus C)$, one needs an
explicit construction of~$C$ (Theorem~\ref{th.6.enum}
merely states the existence); then, an appropriate version of the
Zariski--van Kampen theorem~\cite{vanKampen} can be used. Sextics
with a triple point can be reduced to trigonal curves
(see~\cite{dessin-e7}, \cite{dessin-e8}, \cite{dessin-e6}; the
case of $\bD$ type singularities is pending), and their braid
monodromy, hence fundamental group, can be computed in a purely
combinatorial way (see~\cite{degt.kplets} and
Section~\ref{S.trigonal} below for a further discussion).
Sextics with an involutive symmetry can also be reduced to
trigonal curves (see~\cite{degt.Oka3}, \cite{degt.8a2},
\cite{degt.e6}, \cite{degt.2a8}, \cite{Oka.symmetric}); the
computation becomes slightly more involved as one needs to keep
track of an extra section. One may hope to produce a few more
deformation families using reducible curves and/or unstable
symmetries. For the remaining curves, those with $\bA$ type
singularities only, most promising seems the approach
of~\cite{Artal.Trends}, producing sextics with at least eight
double points, possibly infinitely near. This approach is also
based on a certain involutive symmetry
(see comments to Problem~\ref{prbl.symmetric.br}),
reducing the computation of the
fundamental group
to curves of smaller degree.


The following problem, although not dealing with plane sextics directly,
is closely related to~\ref{prbl.pi1}.

\problem\label{prbl.group.K3}
What are the groups $\pi_1(X\sminus E)$, where $X$ is a
$K3$-surface and $E$ is an $\bA\bD\bE$-configuration, \ie.,
a configuration of $(-2)$-curves spanning a negative definite
sublattice?
\endproblem

If $X=\tX_C$ for a simple sextic $C\subset\Cp2$ and $E$ are the
exceptional divisors contracted in~$\Cp2$, then the group in
question is an index~$2$ subgroup of $\pi_1(\Cp2\sminus C)/\Gd^2$,
where $\Gd$ is the class of a meridian.
Problem~\ref{prbl.group.K3} was posed in~\cite{ShimadaZhang},
where certain conditions for $\pi_1(X\sminus E)$ to be trivial
were found. Numerous examples of nontrivial groups are given by
Xiao's approach~\cite{Xiao} to the classification of finite groups
acting symplectically on $K3$-surfaces.

\subsection{Geometry of simple sextics}\label{s.geometry}
The developed theory of $K3$-surfaces (the global Torelli theorem,
the existence of fine moduli spaces, the Riemann--Roch theorem)
lets one study more subtle
geometric properties of simple sextics, such as symmetries, splitting
curves, minimal fields of definition, \etc. It seems to be a
common belief that, given time, any particular question concerning
$K3$-surfaces can be answered (although Problems~\ref{prbl.real}
and~\ref{prbl.real.diff=def} involving anti-holomorphic
automorphisms seem
more difficult than the purely holomorphic ones). Thus, in my
opinion, the most interesting part is an attempt to understand the
extent to which the properties of (simple) plane sextics
generalize to
curves of higher degree, see \eg. Problem~\ref{prbl.splitting.d}.

\problem\label{prbl.symmetries}
Describe the groups of stable (under equisingular deformations)
symmetries of reducible simple
sextics. Is the quotient of such a sextic by a stable involutive
symmetry always a maximal trigonal/hyperelliptic curve?
\endproblem

By a \emph{symmetry} of a plane curve $C\subset\Cp2$ we mean
a projective transformation of~$\Cp2$ preserving~$C$ as a set.
The answer to the
corresponding question for irreducible simple sextics is given
in~\cite{symmetric}.

If $\mu(C)\le18$, the group of stable symmetries of~$C$ depends on
the configuration $\tS_C\supset S_C$ of~$C$ only: it is the group
of automorphisms of the configuration acting identically on its
discriminant. If $\mu(C)=19$, then $C$ is rigid and all its
symmetries are stable. Hence, the problem can be solved by a
careful analysis of
Shimada's lists of lattice types~\cite{Shimada.MMJ} and
maximizing sextics~\cite{Shimada}.

\problem\label{prbl.symmetric}
Describe the equisingular deformation classes of simple sextics
with a prescribed finite group of symmetries.
\endproblem

\problem\label{prbl.real}
Describe the equivariant
equisingular deformation classes of real simple
sextics or, more generally, those of simple sextics with a
prescribed finite
Klein group of symmetries.
\endproblem

Recall that a \emph{real structure} on an algebraic variety~$X$
is an anti-holomorphic involution $c\:X\to X$. More generally, a
\emph{Klein action} of a group~$G$ on~$X$ is an action of~$G$ by both
holomorphic and anti-holomorphic maps. Thus,
Problem~\ref{prbl.real} is a generalization
of~\ref{prbl.symmetric}.

A characterization of finite groups that may act on $K3$-surfaces
is known, see~\cite{Kondo}, \cite{Kondo2}, \cite{Mukai},
\cite{Nikulin2}, and~\cite{Xiao}, and one can reduce both problems
to the study of equivariant isomorphism classes of
oriented abstract homological types~$(\CH,\go)$ with
a prescribed group~$G$ of (anti-)automorphisms.
(Here, an anti-automorphism is an isometry induced by an
anti-holomorphic map; it
should reverse $h$ and~$\go$ and take~$\Gs$ to $-\Gs$ as a set).
Note, however, that in the presence of anti-holomorphic maps the
corresponding moduli space $\CM_6(\CH,\go;G)$ does not need to be
connected (see~\cite{duke} for an example,
with~$G$ as simple as~$\DG6$; this example is not related to
plane sextics, but it does illustrate the phenomenon).
Hence, an analog of
Theorem~\ref{th.6.enum} does not hold in this case and
Problem~\ref{prbl.real} remains meaningful even in the
following simplified version (and even for nonsingular
sextics).

\problem\label{prbl.real.diff=def}
Given a simple sextic $C\subset\Cp2$ with a finite group Klein
action, is its equivariant equisingular deformation type
determined by the diffeomorphism type of the action on $(\Cp2,C)$?
\endproblem

In the example of~\cite{duke}, the answer to this question is in
the affirmative: the three classes differ by the topology of the
real point set. The only other results concerning real sextics
that I know are the classification of nonsingular
sextics~\cite{Nikulin} and the classification of sextics with a
single node~\cite{Itenberg}.

\problem\label{prbl.symmetric.br}
Describe stable birational symmetries of simple sextics. Find a
lattice theoretic description of such symmetries.
\endproblem

A large number of examples is found in~\cite{Artal.Trends}, where
it is shown that any sextic with at least eight double points,
possibly infinitely near, admits an involutive birational
symmetry. The eight double points are used to define a pencil of
elliptic curves, and the involution acts on the corresponding
rational elliptic surface.

\problem\label{prbl.splitting}
Extend Shimada's classification~\cite{Shimada.MMJ} of
$Z$-splitting curves. For example, are there sextics with pairs,
triples, \etc\. of $Z$-splitting curves? (As an alternative, but
probably less interesting problem, one can
study simple sextics with
unstable splitting curves.)
\endproblem

\problem\label{prbl.splitting.d}
What is the correct generalization of the theory of $Z$-splitting
curves to curves of higher degree? As a wild guess, are
$Z$-splitting curves related to the
(integral) Alexander modules~\cite{Libgober2}
or twisted Alexander polynomials~\cite{Cogolludo}?
\endproblem

Roughly, a nonsingular curve~$B$ is called
\emph{$Z$-splitting} for a simple sextic~$C$ if
\roster
\item
the pull-back~$\tilde B$ of~$B$
in the double covering~$\tX_C$
splits, and
\item
$B$ is stable, \ie., it follows equisingular deformations of~$C$
retaining the splitting property.
\endroster
(If $g(B)>0$, it is required in addition that the two components
of~$\tilde B$ should realize distinct classes in $H_2(\tX_C)$;
this seems to be a purely technical assumption, as well as the
requirement that $B$ should be nonsingular or irreducible.)
The existence of $Z$-splitting curves of small degrees is a
lattice invariant; it is due to the torsion $\tS_C/S_C$, hence to
the existence of certain dihedral coverings (see~\cite{degt.Oka}
or~\cite{Shimada.MMJ}).
Thus, in the case of
simple sextics, $Z$-splitting curves \emph{are} related to the
Alexander module. However, whereas the direct implication
($Z$-splitting curves $\implies$ dihedral coverings) is of a
purely topological nature, the converse relies upon the
Riemann--Roch theorem for $K3$-surfaces and does not generalize
directly to higher degrees. (A similar relation between the
existence of
$Z$-splitting sections and the Alexander module for trigonal
curves is mentioned in Remark~\ref{rem.MW} below.)

One can also speculate that the situation with $Z$-splitting
curves is very similar to Theorem~\ref{th.Delta}, where often the
conclusion that $\Delta_C(t)\ne1$ follows from the fact that one
of the linear systems~$\CL_i$ is nonempty while having negative
virtual dimension, \ie., from the \emph{existence} of certain
auxiliary curves. Furthermore,
these auxiliary curves are stable under deformations, and often
they are also splitting in an appropriate covering (\eg., this is
always the case for simple sextics of torus
type~\cite{Shimada.MMJ}, as well as for many curves considered
in~\cite{Oka.Zariski} and~\cite{Oka.deg10}).

In general, $Z$-splitting curves seem to be
more complicated than those
appearing in Theorem~\ref{th.Delta}: the conditions are not always
linear in terms of the original sextic~$C$. (In some of the
examples found in~\cite{Shimada.MMJ}, $Z$-splitting curves
are tangent to~$C$ at its smooth points.)

Another example illustrating the relation between (unstable)
splitting curves and dihedral coverings was recently discovered
in~\cite{Tokunaga.splitting}, where the r\^ole of plane sextics is
played by trigonal curves in the Hirzebruch
surface~$\Sigma_2$. The splitting sections used
in~\cite{Tokunaga.splitting} are triple tangent to the curve at
its smooth points.

\problem\label{prbl.moduli}
What are the geometric properties of the equisingular moduli
spaces $\CM_6(\CH,\go)$?
\endproblem

In principle, the spaces $\CM_6(\CH,\go)$ are given by
Theorem~\ref{th.moduli}. For example, it follows immediately that
the moduli spaces of maximizing sextics are one point sets. The
two other statements that I know are the following:
\Dashes
\dash
the moduli space of sextics with the set of singularities $3\bE_6$
splits into two rational curves (see~\cite{Oka.torus.2});
\dash
each moduli space
containing an irreducible sextic with a stable symmetry is
unirational (see~\cite{degt.2a8}, \cite{degt.Oka3}, \cite{degt.8a2},
\cite{degt.e6}, and~\cite{degt-Oka} for a case by case analysis).
\endDashes
These results are proved geometrically rather than using
Theorem~\ref{th.moduli} directly.

\problem\label{prbl.Galois}
What are minimal fields of definition of maximizing sextics? Are
maximizing simple sextics sharing the same configuration and
transcendental lattice always Galois conjugate?
\endproblem

In all examples where explicit equations are known (\eg.,
\cite{Artal.braids}, \cite{Artal.Trends}, \cite{Artal.KT},
\cite{EyralOka}, \cite{Shimada.AZP};
see~\cite{Artal.survey} for a more complete
list), the curves \emph{are} Galois conjugate and they can be
defined over an algebraic number field of minimal degree (equal to
the number of conjugate curves). General upper and lower bounds to
the degree of the field of definition of a singular $K3$-surface
are found in~\cite{Schutt} and~\cite{Shimada.K3}; however, the
upper bound seems too weak to provide the minimal degree for simple
sextics.

\subsection{Other problems related to $K3$-surfaces}\label{s.others}
There is another class of plane curves related to $K3$-surfaces:
the ramification loci of a generic projection to~$\Cp2$ of a
$K3$-surface $X\subset\Cp{N}$. Each deformation family is uniquely
determined by a pair $n,k$ of positive integers,
so that the hyperplane section~$h$ has square $2nk^2$ (hence
$N=nk^2+1$) and $h/k$ is primitive in $\Pic X$. (There are a few
exceptions corresponding to the hyperelliptic case,
see~\cite{Saint-Donat} for details.)

\problem\label{prbl.projection}
What is the topology of the ramification locus $C\subset\Cp2$ of a
generic projection $X\to\Cp2$ of a $K3$-surface $X\subset\Cp{N}$?
What is the fundamental group $\pi_1(\Cp2\sminus C)$? The group
$\pi_1(X^G)$ of the Galois closure $X^G\to X\to\Cp2$?
\endproblem

A few results in this direction were recently obtained
in~\cite{TeicherFriedman}, where the ramification loci were
studied by degenerating $K3$-surfaces to unions of planes.
As a common generalization, one can also state the following
problem.

\problem\label{prbl.projection.sing}
What is the topology of the ramification locus $C\subset\Cp2$ of a
generic projection $X\to\Cp2$ of a \emph{singular}
$K3$-surface $X\subset\Cp{N}$?
\endproblem

In what concerns the classification, one can prove a theorem
similar to~\ref{th.6.enum}, replacing~$h$ with $h^2=2$ in
Definition~\ref{def.config} with a primitive (in~$\tS$)
element~$h'$ with $(h')^2=2n$ and adjusting
condition~\iref{def.config}{config.not.elliptic} appropriately,
see~\cite{Saint-Donat}. However, it is not quite clear what
the ramification locus is, especially when the surface is
rigid.

\section{Trigonal curves and elliptic surfaces\label{S.trigonal}}

Let $C\subset\Cp2$ be a plane curve of degree~$d$ with a
distinguished singular
point~$S$ of multiplicity $(d-3)$ (\eg., a sextic with a triple
point). Blow~$S$ up to obtain the ruled surface
$\Cp2(S)\cong\Sigma_1$. Then, the proper transform of~$C$ is a
curve $C'\subset\Sigma_1$ intersecting each fiber of the ruling at
three points, \ie., a generalized trigonal curve.
In this section, we outline a combinatorial approach to
such curves and covering elliptic surfaces (which play the
r\^ole of covering $K3$-surfaces for simple sextics).

For the sake of simplicity, in this paper
we consider the case of {\em rational base only}.
Most results cited below
extend more or less directly to trigonal curves and
elliptic surfaces over an arbitrary base~$B$ provided that it is
considered as a {\em topological surface}, \ie., the analytic
structure on~$B$ is allowed to vary during the deformations.

\subsection{Trigonal curves}\label{s.trigonal}
The \emph{Hirzebruch surface}
is a geometrically ruled
rational surface $\Sigma_k\to\Cp1$
with an \emph{exceptional section}~$E$
of square~$-k$, $k\ge0$.
A \emph{generalized trigonal curve} is a curve
$C\subset\Sigma_k$
intersecting each generic fiber at three points.
A \emph{trigonal curve} is a generalized trigonal curve disjoint
from~$E$.

In what follows, {\em we assume that the curves do not contain
as components fibers of the ruling}.
Furthermore, given a generalized trigonal curve $C\subset\Sigma_k$,
the points
of intersection $C\cap E$ can be removed by a sequence of
elementary transformations $\Sigma_k\dashrightarrow\Sigma_{k+1}$.
Similarly, all
non-simple singular points of~$C$ can be converted to
simple ones by inverse elementary transformations
$\Sigma_k\dashrightarrow\Sigma_{k-1}$.
For this reason, {\em we will consider
simple trigonal curves only}. (In fact, it would even suffice to
deal with trigonal curves with double singular points only.
Certainly,
in the general case
one would need to keep track of the vertical components and the
elementary transformations used.)

If $k$ is even, the minimal resolution of singularities~$\tX_C$ of
the double covering of~$\Sigma_k$ ramified at~$C$ and~$E$ is a
relatively minimal Jacobian elliptic surface. Conversely, given a
Jacobian elliptic surface~$X$, the quotient $X/\pm1$ contracts
to~$\Sigma_k$ and is ramified at~$E$ and a trigonal curve~$C$.
With an abuse of the language, we call~$\tX_C$ the elliptic surface
ramified at~$C$ and~$C$ the ramification locus of~$X$. If $k$ is
odd,
the surface~$\tX_C$ is only
defined locally with respect to the base of the ruling.
A \emph{singular fiber} of~$C$ is the projection of a singular
fiber of~$\tX_C$. For the topological
types of singular fibers, we use the common
notation referring to the extended Dynkin graphs of exceptional
divisors (see Table~\ref{tab.j}), the advantage being the fact
that it reflects the types of the singular points of~$C$.
The fibers of type~$\tA_0^{**}$, $\tA_1^*$, and~$\tA_2^*$ are
called \emph{unstable}; their types are not
necessarily preserved under equisingular
deformations of~$C$.

\midinsert
\table\label{tab.j}
Types of singular fibers
\endtable
\def\tabrule{\vrule height11pt depth4pt}\def\neg{\!\!\!}
\def\I{\roman{I}}\def\II{\roman{II}}\def\III{\roman{III}}\def\IV{\roman{IV}}

\centerline{\vbox{\offinterlineskip
\halign{\tabrule\quad\hss#\hss\quad\vrule&&\quad\hss$#$\hss\quad\vrule\cr
\noalign{\hrule}
\multispan2\tabrule\hss Type of~$F$\hss\vrule&
 \neg j(F)\neg&\neg\text{Vertex}\neg&\neg\text{Valency}\neg\cr
\noalign{\hrule}
$\tilde\bA_p$ ($\tilde\bD_{p+5}$), $p\ge1$&
 \I_{p+1}\ (\I^*_{p+1})&
 \infty&\text{\cross-}&2(p+1)\cr
\noalign{\hrule}
$\tilde\bA_0^*$ ($\tilde\bD_5$)&\I_1\ (\I^*_1)&
 \infty&\text{\cross-}&2\cr
\noalign{\hrule}
$\tilde\bA_0^{**}$ ($\tilde\bE_6$)&\II\ (\II^*)&
 0&\text{\black-}&2\bmod6\cr
\noalign{\hrule}
$\tilde\bA_1^*$ ($\tilde\bE_7$)&\III\ (\III^*)&
 1&\text{\white-}&2\bmod4\cr
\noalign{\hrule}
$\tilde\bA_2^*$ ($\tilde\bE_8$)&\IV\ (\IV^*)&
 0&\text{\black-}&4\bmod6\cr
\noalign{\hrule}\crcr}}}
\eightpoint\botcaption{\bf Comments}
Fibers of type~$\tilde\bA_0$ are not singular;
fibers of type~$\tilde\bD_4$
are
not detected by the $j$-invariant.
Fibers of type~$\tA_0$ or~$\tD_4$
with complex multiplication of order~$2$
(respectively,~$3$) are over the \white-vertices of valency
$0\bmod4$ (respectively,
over the
\black-vertices of valency $0\bmod6$). The types shown
parenthetically
are obtained from the corresponding $\tA$-types
by an elementary transformation.
\endcaption
\endinsert

The \emph{$j$-invariant} $j_C\:\Cp1\to\Cp1$ is the
analytic continuation of the function
sending a nonsingular fiber of~$\tX_C$ to its $j$-invariant (divided
by~$12^3$). By definition, $j_C$ is a birational invariant; hence
it extends to generalized/non-simple trigonal curves. If
$j_C=\const$, the curve is called \emph{isotrivial}; such curves
form a very restricted class. From now on, we assume all curves
non-isotrivial.

\definition\label{def.max}
A non-isotrivial trigonal curve~$C$
is called \emph{maximal} if it has the following properties:
\roster
\item\local{noD4}
$C$ has no singular fibers of type~$\tD_4$;
\item\local{0,1,infty}
$j=j_C$ has no critical values other than~$0$, $1$, and~$\infty$;
\item\local{le3}
each point in the pull-back $j^{-1}(0)$ has ramification index at
most~$3$;
\item\local{le2}
each point in the pull-back $j^{-1}(1)$ has ramification index at
most~$2$.
\endroster
A Jacobian elliptic surface is \emph{maximal}
if its ramification
locus is maximal. A maximal surface is \emph{extremal} if it
has no unstable fibers.
\enddefinition

The more topological definition of extremal elliptic surfaces
given above fits better into the framework of this survey;
in fact, it is the content of~\cite{MNori}.

Maximal trigonal curves and maximal elliptic surfaces
are defined over algebraic number fields.
Maximal curves are indeed maximal in the sense of the
total Milnor number, see~\cite{dessin-e7}: for a
non-isotrivial curve $C\subset\Sigma_k$ one has
$$
\mu(C)\le5k-2-\#\{\text{\rm unstable fibers of~$C$}\},
$$
the equality holding
if and only if $C$ is maximal. (Thus, truly maximal are
curves without unstable fibers, \ie., those corresponding to
extremal elliptic surfaces.)

\subsection{Dessins}\label{s.dessins}
A \emph{trichotomic graph}
is
a
directed
graph $\Gamma\subset\base$
decorated with the
following additional structures
(called the \emph{colorings} of the edges and vertices
of~$\Gamma$):
\Dashes
\dash
each edge of~$\Gamma$ is of one of the three kinds: solid, bold,
or dotted;
\dash
each vertex of~$\Gamma$ is of one of the four kinds: \black-, \white-,
\cross-, or monochrome
(the vertices of the first three kinds being called \emph{essential})
\endDashes
and satisfying the following conditions:
\roster
\item\local{tg-valency}
the valency of each essential vertex
is at least~$2$,
and the valency of each monochrome vertex
is at least~$3$;
\item\local{tg-oriented}
the orientations of the edges of~$\Gamma$ form an orientation
of
$\partial(\base\sminus\Gamma)$;
\item\local{tg-monochrome}
all edges incident to a monochrome vertex are of the same kind;
\item\local{tg-cross}
\cross-vertices are incident to incoming dotted edges and
outgoing solid edges;
\item\local{tg-black}
\black-vertices are incident to incoming solid edges and
outgoing bold edges;
\item\local{tg-white}
\white-vertices are incident to incoming bold edges and
outgoing dotted edges.
\endroster
In \loccit{tg-cross}--\loccit{tg-white} the lists are complete,
\ie., vertices cannot be incident to edges of other kinds or with
different orientation. \black-vertices of valency $0\bmod6$ and
\white-vertices of valency $0\bmod4$ are called
\emph{nonsingular}; all other essential vertices are
\emph{singular}.


A path in a trichotomic graph~$\Gamma$ is
\emph{monochrome} if all its vertices are
monochrome.
The graph~$\Gamma$ is said to be
\emph{admissible} if it has no oriented monochrome cycles.
A \emph{dessin} is an admissible trichotomic graph.

Let $\Gamma\subset\base$ be a trichotomic graph, and let~$v$ be a
vertex of~$\Gamma$. Pick a regular neighborhood $U\ni v$
and replace the intersection $\Gamma\cap U$ with another
decorated graph, so that the result~$\Gamma'$ is again a
trichotomic graph. If $\Gamma'\cap U$ contains essential
vertices of at most one kind
and
contains no monochrome vertices, then $\Gamma'$ is called a
\emph{perturbation} of~$\Gamma$ (at~$v$),
and the original graph~$\Gamma$ is called a
\emph{degeneration} of~$\Gamma'$.
The perturbation~$\Gamma'$ as above is called
\emph{equisingular} if
the intersection $\Gamma'\cap U$ has at most one singular
vertex. (In this case, $\Gamma$ is an
\emph{equisingular degeneration} of~$\Gamma'$.)
Two dessins $\Gamma',\Gamma''\subset\base$
are
said to be \emph{equivalent} if they can be connected by a
chain $\Gamma'=\Gamma_0,\Gamma_1,\ldots,\Gamma_n=\Gamma''$ of
dessins
so that
each~$\Gamma_i$, $1\le i\le n$, either is isotopic
to~$\Gamma_{i-1}$
or is an equisingular perturbation or
degeneration of~$\Gamma_{i-1}$. Clearly, equivalence of dessins is
an equivalence relation.

Following Orevkov~\cite{Orevkov}, define the dessin~$\Gamma_C$ of
a trigonal curve~$C$ as follows.
As a set, $\Gamma_C=j_C^{-1}(\Rp1)$;
the \black--, \white--,
and \cross-vertices are  the pull-backs of~$0$, $1$, and~$\infty$,
respectively (monochrome vertices being the
critical points of~$j_C$
with other real critical values), the edges are solid, bold,
or dotted provided that their images belong to $[\infty,0]$,
$[0,1]$, or $[1,\infty]$, respectively, and the orientation
of the edges is that induced from the positive orientation
of~$\Rp1$ (\ie., order of~$\R$).

The relation between the vertices of~$\Gamma_C$ and singular
fibers of~$C$ is outlined in Table~\ref{tab.j};
it is this
relation that motivates the above
definition of \emph{equisingular}
perturbation of dessins.

Next theorem,
based on the Riemann existence theorem, is motivated
by~\cite{Orevkov}. Its proof
is essentially contained on~\cite{DIK.elliptic}, see
also~\cite{degt.kplets}. A similar statement holds for
generalized trigonal curves and those with triple
or even non-simple singular points:
one merely needs to consider
dessins with appropriately
marked vertices in order to keep track of
the elementary transformations, \cf. Subsection~\ref{s.trigonal}.

\theorem\label{th.tc.enum}
The map $C\mapsto\Gamma_C$
sending a trigonal curve~$C$ to its dessin
establishes a
bijection
between the set of
fiberwise deformation classes of
trigonal curves
\emph{with double singular points only}
and the set of equivalence classes of dessins.
\endtheorem

As a consequence, one obtains the following important property of
trigonal curves that makes them similar to plane sextics (\cf.
Theorem~\ref{th.perturbations}): {\proclaimfont the singular
fibers of a trigonal curve can be perturbed arbitrarily and
independently}.
This statement can be proved by the standard `cut-and-paste'
techniques used to show that any elliptic surface
deforms
to a generic one. Alternatively, it can easily be
proved  using dessins.
Note though that,
when a triple point or an unstable fiber is perturbed, the
$j$-invariant may change discontinuously, increasing its degree,
so that a single vertex of the dessin (the one representing the
fiber perturbed)
is removed and replaced with a fragment
containing essential vertices of all three kinds.

\subsection{Skeletons}\label{s.skeleton}
Dessins of maximal trigonal curves are called \emph{maximal}.
Such dessins do not admit nontrivial degenerations. Therefore, two
maximal dessins are equivalent if and only if they are isotopic.
A convenient way to encode maximal dessins is using
\emph{skeletons}, which are in fact the \emph{dessins d'enfants}
in Grothendieck's original sense. By definition, the
skeleton~$\Sk_C$ is obtained from the dessin~$\Gamma_C$ by
removing all \cross-vertices, solid and dotted edges, and
\white-vertices of valency~$2$. The skeleton $\Sk_C$ of a maximal
trigonal curve~$C$ has the following properties:
\roster
\item\local{Sk.connected}
$\Sk_C$ is connected;
\item\local{Sk.3}
the valency of each \black-vertex is at most~$3$;
\item\local{Sk.1}
the valency of each \white-vertex is~$1$;
\item\local{Sk.black}
there is at least one \black-vertex.
\endroster
Conversely, any graph $\Sk\subset\base$ satisfying
conditions~\loccit{Sk.connected}--\loccit{Sk.black} above
extends to a unique
maximal dessin: place a \white-vertex at the center of each edge
connecting two \black-vertices, place a \cross-vertex $v_R$ at the
center of each region~$R$ of
the complement $\base\sminus\Sk$, and connect $v_R$
to the vertices in $\partial R$ by appropriate edges in the star
like manner. (Note that $R\cong\operatorname{Cone}\partial R$ is
an open
topological disk due to~\loccit{Sk.connected}.)
Thus, next statement is an immediate consequence of
Theorem~\ref{th.tc.enum}.

\theorem\label{th.Sk.enum}
The map $C\mapsto\Sk_C$ establishes a one-to-one correspondence
between the set of fiberwise deformation classes of maximal
trigonal curves with double singular points only and the set of
orientation preserving homeomorphism classes of graphs
$\Sk\subset\base$ satisfying
conditions~\itemref{s.skeleton}{Sk.connected}--\ditto{Sk.black}
above.
\endtheorem

\subsection{Classification of dessins}\label{s.dessins.enum}
As in the case of plane sextics, we start with a few
classification problems, some easy, some hopeless.

\problem\label{prbl.invariants}
Find combinatorial invariants of dessins.
\endproblem

Since the equivalence of dessins defined in
Subsection~\ref{s.equivalence}
involves
perturbations and degenerations,
it is difficult to decide whether two non-maximal
dessins are equivalent.
One approach would be to try to
compute the fundamental groups, Alexander modules/polynomials,
transcendental lattices, and
similar
invariants of the
related trigonal curves and elliptic surfaces, see
Subsections~\ref{s.tc.topology} and~\ref{s.es.topology}.
However,
known
invariants do not always distinguish dessins, even
maximal, see~\ref{s.es.topology}.

\midinsert
\centerline{\picture{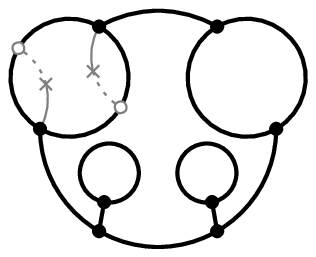}}
\figure
Skeleton~$\Sk$ (black) for Problem~\ref{prbl.example}
\endfigure\label{fig.e6+a11+a1}
\endinsert

As another example, start with the skeleton~$\Sk$ shown in black in
Figure~\ref{fig.e6+a11+a1} and convert it to a dessin~$\Gamma$ by
placing \emph{two} monovalent \cross-vertices inside one of the
two bigons (as shown in grey in Figure~\ref{fig.e6+a11+a1})
and a single \cross-vertex inside each other region (\cf.
Subsection~\ref{s.skeleton}).
Let $C\subset\Sigma_4$ be a trigonal curve whose dessin
is~$\Gamma$. It is birationally equivalent (the transformation
being uniquely determined by the curve) to an irreducible plane
sextic~$C'$ not of torus type and with the set of singularities
$\bE_6\splus\bA_{11}\splus\bA_1$. Such sextics are known to form
two complex conjugate equisingular strata, see~\cite{Aysegul}.

\problem\label{prbl.example}
Find a direct combinatorial proof of the fact that the
dessin~$\Gamma$ described above is not equivalent to its mirror
image.
\endproblem

\problem\label{prbl.skeletons}
Classify maximal dessins.
\endproblem

This problem is too vague and too hopeless. Even the simple
counting of the number of skeletons with given sizes of the
regions is a very difficult task: it is a special case of the so
called \emph{triple Hurwitz numbers}. The study of the Hurwitz
numbers is a separate intensively developing area of mathematics:
they are related to characters of symmetric groups,
moduli spaces of stable curves,
Gromov--Witten invariants, and more. For further
references on Hurwitz numbers, see~\cite{Goulden}
or~\cite{Vainshtein}. Even the asymptotic behavior of these
numbers is of a certain interest: for example,
one can construct exponentially large
Zariski $k$-plets of trigonal curves (see~\cite{degt.kplets} and
comments to Theorem~\ref{th.es.tree} below).

The complete classification of skeletons with $2k$ vertices for
$k\le2$ is easy and well known (see, \eg.,~\cite{symmetric}). A
partial classification for $k=3$ or~$4$ is found
in~\cite{dessin-e7}, \cite{dessin-e8}, and~\cite{dessin-e6}, where
only the dessins arising from sextics with an $\bE$ type singular
point are considered. Next, and probably last, case would be that
of $K3$-surfaces.

\problem\label{prbl.K3.elliptic}
Classify extremal elliptic $K3$-surfaces.
\endproblem

This is equivalent to the study of maximal trigonal curves
$C\subset\Sigma_k$ with $t$ triple singular points such that
$k+t=4$. The case of $t=0$
and stable singular fibers only
is contained in~\cite{Beukers}, and the
case $t\ge2$ follows from the known results for $k\le2$. The
remaining case $k=3$, $t=1$ is still pending.

\subsection{Topology of trigonal curves}\label{s.tc.topology}
The
class realized by a trigonal curve $C\subset\Sigma_k$
is $[C]=3[E]+3k[F]$, where $E$ and~$F$ are, respectively, the
exceptional section and a generic fiber of~$\Sigma_k$. Hence,
there is a unique cyclic covering of~$\Sigma_k$ ramified at~$C$,
and this covering is $3$-fold. Although the latter might be well
worth studying, it is nevertheless customary, when speaking about
the fundamental group, Alexander module, \PL-equivalence,
\etc\. of a (generalized)
trigonal curve $C\subset\Sigma_k$, to {\em refer to the corresponding
invariants of the complement $\Sigma_k\sminus(C\cup E)$ or
pair $(\Sigma_k,C\cup E)$}.

\problem\label{prbl.tc.pi1}
What finite groups~$G$ can appear as the fundamental groups of
(generalized) trigonal curves? What finite commutants $[G,G]$ can
appear?
\endproblem

The corresponding problem for hyperelliptic curves in Hirzebruch
surfaces (hence for plane curves of degree~$d$ with a singular
point of multiplicity $(d-2)$\,) has been solved in~\cite{groups}.
Apart from the cyclic groups~$\CG{r}$, only three finite groups
can appear as the commutant $[G,G]$:
$\QG$,
$\SL(2,\FF_3)$, and $\SL(2,\FF_5)$.

In the case of trigonal curves, the braid monodromy, hence a
presentation of~$G$, are determined by the dessin in a simple
combinatorial way, see~\cite{degt.kplets}. (In the case of
generalized trigonal curves, additional care should be taken about
the points at the exceptional section, \cf.~\cite{dessin-e7}.)
However, the situation
is more complicated as, first, there is a huge number of
non-equivalent dessins (see comments to
Problem~\ref{prbl.skeletons}) and, second, the monodromy takes
values in the non-abelian group~$\BG3$, which makes the
further analysis
of the presentations
rather difficult.
As far as I know, only very few special cases (mainly those
arising from plane sextics, see~\cite{dessin-e7}, \cite{dessin-e8},
and~\cite{dessin-e6})
have been tried so far. The finite commutants discovered are
$\CG5$, $\SL(2,\FF_5)$, and $\SL(2,\FF_{19})$. (For reducible
curves, one can also encounter $\CG3$ and $Q_8$.)
Among other examples
is a generalized trigonal curve in~$\Sigma_1$
that blows down to an
irreducible quintic with the set of singularities $3\bA_4$;
the commutant of its fundamental
group is a non-split central extension of $(\CG2)^2$
by~$(\CG2)^4$, see~\cite{groups}.

An answer to Problem~\ref{prbl.tc.pi1}
would shed a new light on the general
problem of realizability of finite groups as fundamental groups of
plane curves, see~\cite{Libgober.survey.1},
\cite{Libgober.survey.2}.
Certainly, this approach has its limitations; for example, all
groups admit a presentation with at most three generators.
(Note though that the natural basis for~$G$
given by the Zariski--van Kampen method
using the ruling of~$\Sigma_k$
is not necessarily
a $C$-basis in the sense of
Kulikov~\cite{Kulikov}.)

\problem\label{prbl.tc.Delta}
Compute the
Alexander like invariants, such as
the Alexander polynomial \cite{Libgober1},
Alexander module~\cite{Libgober2},
other invariants based on representations of the braid
group~\cite{Libgober.Burau},
\etc\.
of a generalized trigonal curve in terms of its dessin.
\endproblem

Of course, one can compute the fundamental group (or even braid
monodromy) and use Fox calculus or other appropriate
tools of combinatorial group theory. However, for simple
invariants, one would expect a simple combinatorial expression.

It is worth mentioning
that there is an analog of Theorem~\ref{th.Delta} for trigonal
curves (as well as for curves on any surface). At present, it
is not quite clear how the linear systems~$\CL_i$ fit into the
description of curves in terms of dessins.

\problem\label{prbl.CL}
Find a relation between the linear systems~$\CL_i$ and dessins.
\endproblem

\problem\label{prbl.tc.Zariski}
Do
maximal trigonal curves sharing the same combinatorial type of
singular fibers but not related by a fiberwise equisingular
deformation
differ topologically,
\ie., do they form Zariski pairs in the sense of
Artal~\cite{Artal.def}?
\endproblem

As mentioned above, there are large collections of
not
deformation equivalent maximal trigonal curves sharing the same
combinatorial type. A few examples, growing exponentially in the
degree, are found in~\cite{degt.kplets}.
In Subsection~\ref{s.es.topology} below it is explained that these
curves share the same fundamental group and homological type; I do
not know whether they differ topologically.

\subsection{Topology and geometry of elliptic surfaces}\label{s.es.topology}
We conclude with a few problems related to the elliptic surfaces
ramified at trigonal curves (which can be regarded as analogues of
the $K3$-surfaces ramified at simple sextics). In addition to a
few general questions, we announce a new result
(Theorem~\ref{th.es.tree}) and discuss related problems concerning
a special class of extremal elliptic surfaces, namely
those defined by
the so called pseudo-trees.

\problem\label{prbl.es.Zariski}
Do
extremal elliptic surfaces
sharing the same combinatorial type of singular fibers but not
related by a fiberwise equisingular deformation
differ topologically?
\endproblem

This problem is closely related to~\ref{prbl.tc.Zariski}: one can
try to use double covering to distinguish curves. In analogy to
simple sextics, define the \emph{homological type} of a
relatively minimal Jacobian
elliptic surface~$X$ as the pair of lattice extensions
$$
H_2(X)\supset\tS_X\supset S_X:=\Z\Gs_X+\Z[E_X]+\Z[F_X],
$$
where
\Dashes
\dash
$\Gs_X$ is the set of classes realized by the $(-2)$-curves in the
fibers of~$X$,
\dash
$E_X$ and~$F_X$ are, respectively, the section and a generic fiber
of~$X$, and
\dash
$\tS_X=H_2(X)\cap(S_X\otimes\Q)$ is the primitive hull of~$S_X$.
\endDashes
(It is no longer required that each root of
$\tS_X\cap\Q\Gs_X$ should belong to $\Z\Gs_X$.)
An \emph{isomorphism} of homological types
is
a bijective isometry of~$H_2(X)$ (regarded as an abstract lattice)
respecting~$\Gs_X$ (as a set), $[E_X]$,
and~$[F_X]$.

\Remark\label{rem.MW}
Another similarity between trigonal curves and plane sextics is a
relation between the existence of splitting ($Z$-splitting?)
sections and the Alexander module, \cf. comments to
Problem~\ref{prbl.splitting.d}. (Here, by a \emph{section} we mean
a section of the Hirzebruch surface~$\Sigma_k$ containing the
curve.)
Indeed, according to~\cite{Shioda},
the quotient $\tS_X/S_X$ equals the torsion of
the Mordell--Weil group $\MW{X}$, and each nontrivial
element of~$\MW{X}$,
unless it is a component of the ramification locus~$C$ of~$X$,
projects
to a splitting section of~$C$. On the other hand,
a simple homological computation shows that
the same quotient $\tS_X/S_X$ controls
the existence of certain dihedral coverings ramified at $C+E$,
\cf.~\cite{degt.Oka}, \cite{degt.Oka2}, or~\cite{Shimada.MMJ}.
\endRemark

The orthogonal complement $T_X:=\tS_X^\perp$ is called the
\emph{transcendental lattice} of~$X$;
it is positive definite if and only if $X$ is extremal.
The subset $\Gs'_X:=\Gs_X\cap[E_X]^\perp$ is a root system; it
encodes the combinatorial type of singular fibers of~$X$
(assuming that $X$ has no unstable fibers).
According to Nikulin's theory of discriminant
forms~\cite{Nikulin}, the isomorphism class of a homological type
is determined by~$\Gs'_X$, $T_X$ (which should `match'~$\Gs'_X$),
and a certain set of finite data that determines the finite index
lattice extension $H_2(X)\supset S_X\oplus T_X$.

\problem\label{prbl.HT}
Describe the homological types of extremal elliptic surfaces; in
particular, describe their transcendental lattices.
\endproblem

A few attempts to attack this problem were recently made
in~\cite{ArimaShimada} and~\cite{degt.TL}. Some of the results
of~\cite{degt.TL} are rather discouraging: one has
the following theorem.

\theorem\label{th.es.tree}
Let~$X$ be a Jacobian elliptic surface with
one of the following combinatorial types of singular fibers\rom:
\roster
\item\local{0.even}
$\tA_{10s-2}\oplus(2s+1)\tA_0^*$, $s\ge1$\rom;
\item\local{0.odd}
$\tD_{10s-2}\oplus(2s)\tA_0^*$, $s\ge1$\rom;
\item\local{1.even}
$\tD_{10s+3}\oplus\tD_5\oplus(2s)\tA_0^*$, $s\ge1$\rom;
\item\local{1.odd}
$\tA_{10s-7}\oplus\tD_5\oplus(2s-1)\tA_0^*$, $s\ge1$.
\endroster
Then, within each of the four series, the isomorphism class of the
transcendental lattice $T_X$ is determined by~$s$ only.
\endtheorem

In Theorem~\ref{th.es.tree}, one
has $T_X\cong-\bD_{2s-2}$ in case~\itemref{th.es.tree}{0.odd}
and $T_X\cong-\bD_{2s-1}\oplus\Z w$, $w^2=4$, in
case~\ditto{1.even} (where, as usual, we let
$\bD_0=0$, $\bD_1=[-4]$, $\bD_2=2\bA_1$, and
$\bD_3=\bA_3$\rom). A precise description of~$T_X$ in the other
two cases is also known, but it is more complicated.

The surfaces in Theorem~\ref{th.es.tree} are extremal. The
skeleton of each surface is a \emph{pseudo-tree}, \ie., it is
obtained from a plane tree $\tree\subset\base$ with all vertices of
valency~$3$ (\emph{nodes}) or~$1$ (\emph{leaves}) by patching each
leaf with a small loop (see Figure~\ref{fig.tree}, where the
original tree~$\tree$
is shown in black and the loops attached, in grey).
The four series in the theorem differ by the number of vertices
in~$\tree$
($4s$ in cases~\ditto{0.even} and~\ditto{1.even}
or $4s-2$ in cases~\ditto{0.odd} and~\ditto{1.odd}\,) and the
number of type~$\tD_5$ singular fibers over the \cross-vertices
in the regions inside the loops attached
to the tree
(\cf. the passage from a skeleton to a dessin in
Subsection~\ref{s.skeleton}; the number in question is
none in cases~\ditto{0.even} and~\ditto{0.odd} or one in
cases~\ditto{1.even} and~\ditto{1.odd}\,).

\midinsert
\centerline{\picture{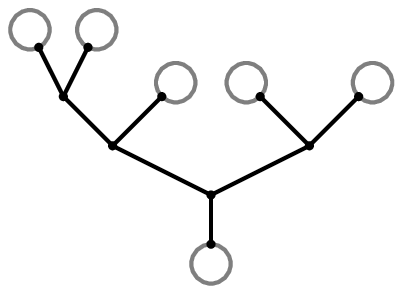}}
\figure
Converting a tree~$\tree$ (black) to a skeleton
\endfigure\label{fig.tree}
\endinsert

The number of
pseudo-trees
grows exponentially in~$s$;
for example, in cases \iref{th.es.tree}{1.even}
and~\ditto{1.odd}, where the possible automorphisms of the tree
are eliminated by distinguishing one of the leaves,
the count is given by the Catalan numbers
$C(2s-1)$ and $C(2s-2)$, respectively,
see~\cite{degt.kplets} for details. On the other hand, the number
of automorphisms of the discriminant form $\discr T_X$
(which may lead
to distinct homological types) grows at most linearly. (This
number is given by $2^m$, where
$m\le\log\,\mathopen|\discr T_X\mathclose|$
is the number of distinct primes dividing
$\mathopen|\discr T_X\mathclose|$.)
Thus, we
obtain exponentially large collections of extremal elliptic
surfaces sharing the same homological type.

We conclude with three problems specific to pseudo-trees.

\problem\label{prbl.tree.D+}
What are the transcendental lattices of extremal elliptic surfaces
defined by pseudo-trees and with more than one singular fiber of
type~$\tD_5$? Are they determined by the number of vertices and
$\tD_5$ type fibers only?
\endproblem

\problem\label{prbl.tree.Galois}
Are the surfaces in Theorem~\ref{th.es.tree} Galois conjugate
(within a fixed series and for a fixed value of~$s$)?
\endproblem

\problem\label{prbl.tree.braid}
Are the braid monodromies of maximal trigonal curves defined by
pseudo-trees Hurwitz equivalent (see \eg.~\cite{Kulikov})?
\endproblem

It is shown in~\cite{degt.TL} that, with very few exceptions,
the fundamental groups of maximal trigonal curves defined by
pseudo-trees are abelian, hence they do not distinguish the braid
monodromies.

For a simple trigonal curve $C\subset\Sigma_k$,
the Hurwitz equivalence class of its
braid monodromy determines and is determined by the fiberwise
diffeomorphism class of the pair $(\Sigma_k,C\cup E)$,
\cf.~\cite{Carmona} and~\cite{KulikovTeicher}.
In general, the problem of distinguishing the Hurwitz equivalence
classes seems very difficult, but in this particular case, where
the
monodromy takes values in the relatively simple
group~$\BG3$, one may hope to get a reasonable
solution.
For example, the transcendental lattice~$T_X$ defined above
can be reconstructed
directly from the braid monodromy, \cf.~\cite{ArimaShimada}
and~\cite{degt.TL},
and it is obviously invariant under Hurwitz equivalence. Thus, one
can generalize Problem~\ref{prbl.tree.braid} as follows.

\problem\label{prbl.braid}
Given a braid monodromy
$\Gb_1,\ldots,\Gb_r\in\BG3$ with all~$\Gb_i$ conjugate to~$\Gs_1$,
is its Hurwitz equivalence class determined by the product
$\Gb_1\ldots\Gb_r\in\BG3$ and the transcendental lattice~$T$?
By the product $\Gb_1\ldots\Gb_r$ only?
\endproblem

In Orevkov's terminology, the last question can be
restated as follows: is the natural map from the $\BG3$
valued braid monodromy monoid to $\BG3$ injective?

\refstyle{C}

\enddocument